# TESTING FOR A LINEAR MA MODEL AGAINST THRESHOLD MA MODELS

By Shiqing Ling[1] and Howell Tong[2]

*Hong Kong University of Science and Technology and London School of Economics and Political Science*

This paper investigates the (conditional) quasi-likelihood ratio test for the threshold in MA models. Under the hypothesis of no threshold, it is shown that the test statistic converges weakly to a function of the centred Gaussian process. Under local alternatives, it is shown that this test has nontrivial asymptotic power. The results are based on a new weak convergence of a linear marked empirical process, which is independently of interest. This paper also gives an invertible expansion of the threshold MA models.

**1. Introduction.** Since Tong [30], threshold autoregressive (TAR) models have become a standard class of nonlinear time series models. Some fundamental results on the probabilistic structure of this class were given by Chan, Petruccelli, Tong and Woolford [11], Chan and Tong [12] and Tong [31]. The 1990s saw many more contributions including, for example, Chen and Tsay [15], Brockwell, Liu and Tweedie [6], Liu and Susko [27], An and Huang [3], An and Chen [1], Liu, Li and Li [26], Ling [23] and others.

The likelihood ratio (LR) test for the threshold in AR models was studied by Chan [8, 9] and Chan and Tong [13]. Tsay [33, 34] proposed some methods for testing the threshold in AR and multivariate models. Lagrange multiplier tests were studied by Wong and Li [35, 36] for (double) TAR–ARCH models. The Wald test was studied by Hansen [17] for TAR models. Testing the threshold in nonstationary AR models was investigated by Caner and Hansen [7]. The asymptotic theory on the estimated threshold parameter in

Received February 2004; revised February 2005.

[1] Supported in part by Hong Kong Research Grants Commission Grants HKUST6113/02P and HKUST4765/03H.

[2] Supported in part by Hong Kong Research Grants Commission Grant HKU7049/03P and by University of Hong Kong Distinguished Research Achievement Award.

TAR models was established by Chan [10] and Chan and Tsay [14]. Recently, Chan's result was extended to non-Gaussian error TAR models by Qian [28]; see also [20] for threshold regression models. Hansen [18] obtained a new limiting distribution for TAR models with changing parameters; see also [19].

However, almost all the research in this area to date has been limited to the AR or AR-type models. Except for Brockwell, Liu and Tweedie [6], Liu and Susko [27], de Gooijer [16] and Ling [23], it seems that threshold moving average (TMA) models have not attracted much attention in the literature. It is well known that, in the linear case, MA models are as important as the AR models. In particular, for many economic data, such as monthly exchange rates, IBM stock market prices and weekly spot rates of the British pound, the models selected in the literature are often MA or ARMA models from the point of view of parsimony; see, for example, [32]. Now, the concept of threshold has been recognized as an important idea for time series modeling. Therefore, it is natural to introduce this concept in the context of MA modeling leading to the TMA models. Again, model parsimony is often an important consideration in nonlinear time series modeling. We shall give an example of this in Section 4. In addition, techniques developed for TMA models should prepare us for a systematic study of the much more challenging threshold ARMA models. We shall give one such instance in the Appendix.

We investigate the quasi-LR test for threshold in MA models. Under the hypothesis of no threshold, it is shown that the test statistic converges weakly to a function of a centred Gaussian process. Under local alternatives, it is shown that this test has nontrivial asymptotic power. The results heavily depend on a linear marked empirical process. This type of empirical process has been found to be very useful and was investigated by An and Cheng [2], Chan [10], Stute [29], Koul and Stute [22], Hansen [18] and Ling [24] for various purposes. However, all the processes in these papers have only one marker. To the best of our knowledge, our linear marked empirical process which includes infinitely many markers has never appeared in the statistical literature before. This is of independent interest. This paper also gives an invertible expansion of the TMA models.

This paper proceeds as follows. Section 2 gives the quasi-LR test and its null asymptotic distribution. Section 3 studies the asymptotic power under local alternatives. Some simulation results and one real example are given in Section 4. Sections 5 and 6 present the proofs of the results stated in Section 2.

## 2. Quasi-LR test and its asymptotics.

The time series $\{y_t : t = 0, \pm 1, \pm 2, \ldots\}$ is said to follow a TMA$(p, q, d)$ model if it satisfies the equation

$$(2.1) \qquad y_t = \sum_{i=1}^{p} \phi_i \varepsilon_{t-i} + \sum_{i=1}^{q} \psi_i I(y_{t-d} \leq r) \varepsilon_{t-i} + \varepsilon_t,$$

where $\{\varepsilon_t\}$ is a sequence of independent and identically distributed (i.i.d.) random variables with mean zero and variance $0 < \sigma^2 < \infty$, $p, q, d$ are known positive integers with $p \geq q$, $I$ is the indicator function and $r \in R$ is called the threshold parameter. Let $\Theta$ and $\Theta_\psi$ be compact subsets of $R^p$ and $R^q$, respectively, and $\Theta_1 = \Theta \times \Theta_\psi$ be the parameter space. Let $\phi = (\phi_1, \ldots, \phi_p)'$, $\psi = (\psi_1, \ldots, \psi_q)'$ and $\lambda = (\phi', \psi')'$. Here $\lambda$ is the unknown parameter (vector) and its true value is $\lambda_0 = (\phi_0', \psi_0')'$. Assume $\lambda_0$ is an interior point in $\Theta_1$.

Given observations $y_1, \ldots, y_n$ from model (2.1), we consider the hypotheses

$$H_0 : \psi_0 = 0 \quad \text{versus} \quad H_1 : \psi_0 \neq 0 \text{ and some } r \in R.$$

Under $H_0$, the true model (2.1) reduces to the usual linear MA model and $\{y_t\}$ is always strictly stationary and ergodic. In this case, the parameter $r$ is absent, which renders the problem nonstandard. Under $H_1$, Liu and Susku [27] and Ling [23] showed that there is always a strictly stationary solution $\{y_t\}$ to the model (2.1) without any restriction on $\lambda_0$. Under $H_0$ and $H_1$, the corresponding quasi-log-likelihood functions based on $\{y_n, y_{n-1}, \ldots\}$ are, respectively,

$$L_{0n}(\phi) = \sum_{t=1}^{n} \varepsilon_t^2(\phi) \quad \text{and} \quad L_{1n}(\lambda, r) = \sum_{t=1}^{n} \varepsilon_t^2(\lambda, r),$$

where $\varepsilon_t(\phi) = \varepsilon_t(\lambda, -\infty)$ and

$$\varepsilon_t(\lambda, r) = y_t - \sum_{i=1}^{p} \phi_i \varepsilon_{t-i}(\lambda, r) - \sum_{i=1}^{q} \psi_i I(y_{t-d} \leq r) \varepsilon_{t-i}(\lambda, r),$$

which is the residual from the TMA model. To make it meaningful, we need to study the invertibility of this model. Assumption 2.1 below is a condition for this.

ASSUMPTION 2.1. $\sum_{i=1}^{p} |\phi_i| < 1$ and $\sum_{i=1}^{p} |\phi_i + \psi_i| < 1$, where $\psi_i = 0$ for $i > q$.

This assumption is similar to Lemma 3.1 for the ergodicity of TAR models in [12]. We discuss the invertibility of a general TMA model in the Appendix.

Since there are only $n$ observations, we need the initial values $y_s$, when $s \leq 0$, to calculate $\varepsilon_t(\phi)$ and $\varepsilon_t(\lambda, r)$. For simplicity, we assume $y_s = 0$ for $s \leq 0$. We denote $\varepsilon_t(\phi)$ and $\varepsilon_t(\lambda, r)$, calculated with these initial values by $\tilde{\varepsilon}_t(\phi)$ and $\tilde{\varepsilon}_t(\lambda, r)$, and modify the corresponding quasi-log-likelihood functions, respectively, to

$$\tilde{L}_{0n}(\phi) = \sum_{t=1}^{n} \tilde{\varepsilon}_t^2(\phi) \quad \text{and} \quad \tilde{L}_{1n}(\lambda, r) = \sum_{t=1}^{n} \tilde{\varepsilon}_t^2(\lambda, r).$$

Let $\tilde{\phi}_n = \arg\min_\Theta \widetilde{L}_{0n}(\phi)$ and $\tilde{\lambda}_n(r) = \arg\min_{\Theta_1} \widetilde{L}_{1n}(\lambda, r)$. We call $\tilde{\phi}_n$ and $\tilde{\lambda}_n(r)$ the conditional least squares estimators of $\phi_0$ and $\lambda_0$, respectively. Given $r$, the quasi-LR test statistic for $H_0$ against $H_1$ is defined as

$$\widetilde{LR}_n(r) = -2[\widetilde{L}_{1n}(\tilde{\lambda}_n(r), r) - \widetilde{L}_{0n}(\tilde{\phi}_n)].$$

Since the threshold parameter $r$ is unknown, a natural test statistic is $\sup_{r \in R} \widetilde{LR}_n(r)$. However, this test statistic diverges to infinity in probability; see (2.2) below and [4]. We consider the supremum of $\widetilde{LR}_n(r)$ on the finite interval $[a, b]$,

$$LR_n = \frac{1}{\tilde{\sigma}_n^2} \sup_{r \in [a,b]} \widetilde{LR}_n(r),$$

where $\tilde{\sigma}_n^2 = \widetilde{L}_{0n}(\tilde{\phi}_n)/n$. This method is used by Chan [8] and Chan and Tong [13]. The idea here is similar to the problem of testing change points in Andrews [4], which has been commonly used in the literature. To study its asymptotics, we need another assumption which is a mild technical condition.

ASSUMPTION 2.2.   $\varepsilon_t$ has a continuous and positive density on $R$ and $E\varepsilon_t^4 < \infty$.

We further introduce the following notation:

$$U_{1t}(\lambda, r) = \partial \varepsilon_t(\lambda, r)/\partial \phi, \qquad U_{2t}(\lambda, r) = \partial \varepsilon_t(\lambda, r)/\partial \psi,$$

$$D_{1t}(\lambda, r) = U_{1t}(\lambda, r)\varepsilon_t(\lambda, r), \qquad D_{2t}(\lambda, r) = U_{2t}(\lambda, r)\varepsilon_t(\lambda, r),$$

$$U_t(\lambda, r) = [U'_{1t}(\lambda, r), U'_{2t}(\lambda, r)]' \quad \text{and} \quad D_t(\lambda, r) = [D'_{1t}(\lambda, r), D'_{2t}(\lambda, r)]'.$$

Throughout this paper, all the expectations are computed under $H_0$. We denote $\Sigma_{rs} = E[U_{2t}(\lambda_0, r)U'_{2t}(\lambda_0, s)]$, $\Sigma_{1r} = E[U_{1t}(\lambda_0, r)U'_{2t}(\lambda_0, r)]$ and $\Omega_r = E[U_t(\lambda_0, r)U'_t(\lambda_0, r)]$. Let $\Sigma = E\{[\partial \varepsilon_t(\phi_0)/\partial \phi][\partial \varepsilon_t(\phi_0)/\partial \phi]'\}$. Here and in the sequel, $o_p(1)$ denotes convergence to zero in probability as $n \to \infty$. We first state one basic lemma, which gives a uniform expansion of $\widetilde{LR}_n(r)$ on $[a, b]$.

LEMMA 2.1.   *If Assumptions* 2.1 *and* 2.2 *hold, then under* $H_0$ *it follows that:*

(a) $\displaystyle \sup_{r \in [a,b]} \|\tilde{\lambda}_n(r) - \lambda_0\| = o_p(1)$,

(b) $\displaystyle \sup_{r \in [a,b]} \left\| \sqrt{n}[\tilde{\lambda}_n(r) - \lambda_0] + \frac{\Omega_r^{-1}}{\sqrt{n}} \sum_{t=1}^n D_t(\lambda_0, r) \right\| = o_p(1)$,

(c) $\displaystyle \sup_{r \in [a,b]} \|\widetilde{LR}_n(r) - T'_n(r)[\Sigma_{rr} - \Sigma'_{1r}\Sigma^{-1}\Sigma_{1r}]^{-1} T_n(r)\| = o_p(1)$,

*where* $T_n(r) = n^{-1/2} \sum_{t=1}^n [D_{2t}(\lambda_0, r) - \Sigma'_{1r}\Sigma^{-1}D_{1t}(\lambda_0, r)]$.

The proof of this lemma is quite complicated and is given in Section 6. Under $H_0$, $D_{1t}(\lambda_0, r) = \varepsilon_t \, \partial \varepsilon_t(\phi_0)/\partial \phi$ and, by (6.4), $D_{2t}(\lambda_0, r)$ has the expansion

$$D_{2t}(\lambda_0, r) = \left[ \sum_{i=0}^{\infty} u' \Phi^i u Z_{t-i-1} I(y_{t-d-i} \leq r) \right] \varepsilon_t \qquad \text{a.s.,}$$

where $Z_t = (\varepsilon_t, \ldots, \varepsilon_{t-q+1})'$, $u = (1, 0, \ldots, 0)'_{p \times 1}$ and $\Phi$ is defined as in Theorem A.1. Following Stute [29], we call $\{T_n(r) : r \in R\}$ a marked empirical process, where each $y_{t-d-i}$ is a marker. It is a linear marked empirical process and includes infinitely many markers. As stated in Section 1, this is a new empirical process. Let $D^q[R_\gamma] = D[R_\gamma] \times \cdots \times D[R_\gamma]$ ($q$ factors), which is equipped with the corresponding product Skorohod topology and in which $R_\gamma = [-\gamma, \gamma]$. Weak convergence on $D^q[R]$ is defined as that on $D^q[R_\gamma]$ for each $\gamma \in (0, \infty)$ as $n \to \infty$ and is denoted by $\Longrightarrow$. We now give the weak convergence of $\{T_n(r) : r \in R\}$ as follows.

THEOREM 2.1. *If Assumption 2.2 holds and all the roots of $z^p - \sum_{i=1}^{p} \phi_i \times z^{p-i} = 0$ lie inside the unit circle, then under $H_0$ it follows that*

$$T_n(r) \Longrightarrow \sigma G_q(r) \qquad \text{in } D^q[R],$$

*where $\{G_q(r) : r \in R\}$ is a $q \times 1$ vector Gaussian process with mean zero and covariance kernel $K_{rs} = \Sigma_{rs} - \Sigma'_{1r} \Sigma^{-1} \Sigma_{1s}$, and almost all its paths are continuous.*

Unlike Koul and Stute [22], our weak convergence does not include the two end-points $\pm \infty$ and $LR_n$ only requires the weak convergence on $D^q[R]$. In addition, our technique heavily depends on $R_\gamma$ and Assumption 2.2. The covariance kernel $K_{rs}$ is essentially different from those of the empirical processes with one marker. Theorem 2.1 is a new weak convergence result and its proof is given in Section 5.

Under $H_0$, it is well known that $\tilde{\sigma}_n^2 = \sigma^2 + o_p(1)$. By Lemma 2.1(c), Theorem 2.1 and the continuous mapping theorem, we obtain the main result as follows.

THEOREM 2.2. *If Assumptions 2.1 and 2.2 hold, then under $H_0$ it follows that*

$$LR_n \overset{\mathcal{L}}{\longrightarrow} \sup_{r \in [a,b]} [G'_q(r) K_{rr}^{-1} G_q(r)]$$

*as $n \to \infty$, where $\overset{\mathcal{L}}{\longrightarrow}$ stands for convergence in distribution.*

When $p = q < d$, $\Sigma_{rr} = \Sigma_{1r} = \Sigma F_y(r)$ since $Z_{t-1}$ and $y_{t-d}$ are independent. Here $F_y(r) = P(y_t \leq r)$. Thus, the limiting distribution is the same as that of

$$(2.2) \qquad \sup_{\beta_1 \leq s \leq \beta_2} \frac{\|B_p(s)\|^2}{s - s^2},$$

where $\beta_1 = F_y(a)$, $\beta_2 = F_y(b)$ and $B_p(s)$ is a $p \times 1$ vector Gaussian process with mean zero and covariance kernel $(r \wedge s - rs)I_p$, where $I_p$ is a $p \times p$ identity matrix. It is interesting that this distribution is the same as that of test statistics for change-points in [4]. The critical values can be found in [4]. In practice, we can select, for example, $\beta_1 = 0.05$ and $\beta_2 = 0.95$. Some guidelines on this can be found in [8]. For given $\beta_1$ and $\beta_2$, we can compute $LR_n$ with $a = F_{ny}^{-1}(\beta_1)$ and $b = F_{ny}^{-1}(\beta_2)$, where $F_{ny}^{-1}(\tau)$ is the $\tau$-quantile of the empirical distribution based on data $\{y_1, \ldots, y_n\}$. For other cases, the critical values of $LR_n$ can be obtained via a simulation method. The implementation is not so difficult in practice.

## 3. Asymptotic power under local alternatives.
To investigate asymptotically the local power of $LR_n$, consider the local alternative hypothesis

$$H_{1n} : \psi_0 = \frac{h}{\sqrt{n}} \qquad \text{for a constant vector } h \in R^q \text{ and } r = r_0 \in R,$$

where $r_0$ is a fixed value. For this, we need some basic concepts as follows. Let $\mathcal{F}^Z$ be the Borel $\sigma$-field on $\mathcal{R}^Z$ with $Z = \{0, \pm 1, \pm 2, \ldots\}$ and $P$ be a probability measure on $(\mathcal{R}^Z, \mathcal{F}^Z)$. Let $P_\lambda^n$ be the restriction of $P$ on $\mathcal{F}_n$, the $\sigma$-field generated by $\{Y_0, y_1, \ldots, y_n\}$, where $Y_0 = \{y_0, y_{-1}, \ldots\}$. Suppose the errors $\{\varepsilon_1(\lambda, r_0), \varepsilon_2(\lambda, r_0), \ldots\}$ under $P_\lambda^n$ are i.i.d. with density $f$ and are independent of $Y_0$. From model (2.1), the distribution of initial value $Y_0$ is the same under both $P_\lambda^n$ and $P_{\lambda_0}^n$. Thus, the log-likelihood ratio $\Lambda_n(\lambda_1, \lambda_2)$ of $P_{\lambda_2}^n$ to $P_{\lambda_1}^n$ is

$$\Lambda_n(\lambda_1, \lambda_2) = 2 \sum_{t=1}^n [\log s_t(\lambda_2) - \log s_t(\lambda_1)],$$

where $s_t(\lambda) = \sqrt{f(\varepsilon_t(\lambda, r_0))}$; see [21] and [25] for details. We first introduce the following assumption.

ASSUMPTION 3.1.   The density $f$ of $\varepsilon_t$ is absolutely continuous with a.e.-derivative and finite Fisher information, $0 < I(f) = \int_{-\infty}^\infty [f'(x)/f(x)]^2 \times f(x)\, dx < \infty$.

The following theorem gives the LAN of $\Lambda_n(\lambda_1, \lambda_2)$ for model (2.1) and the contiguity of $P_{\lambda_0}^n$ and $P_{\lambda_0 + u_n/\sqrt{n}}^n$, where $u_n$ is a bounded constant sequence in $R^{p+q}$.

**Theorem 3.1.** *If Assumptions* 2.1, 2.2 *and* 3.1 *hold and* $\lambda_0 = (\phi_0', 0)'$, *then:*

(a) $\Lambda_n(\lambda_0, \lambda_0 + \frac{u_n}{\sqrt{n}}) = n^{-1/2} u_n' \sum_{t=1}^n U_t(\lambda_0, r_0) \xi_t - I(f) u_n' \Omega_{r_0} u_n / 2 + o_p(1)$ *under* $P_{\lambda_0}^n$, *and*

(b) $P_{\lambda_0}^n$ *and* $P_{\lambda_0 + u_n/\sqrt{n}}^n$ *are contiguous,*

*where* $\xi_t = f'(\varepsilon_t(\lambda_0, r_0))/f(\varepsilon_t(\lambda_0, r_0))$ *and* $\Omega_r$ *is defined as in Lemma* 2.1.

PROOF. By verifying the conditions in Theorem 2.1 and (2.2) in [25], we can show that the conclusions hold. The details are omitted. □

Using Theorem 2.1 and following a routine argument, we can obtain the following theorem. This theorem shows that $LR_n$ has nontrivial local power under $H_{1n}$.

**Theorem 3.2.** *If Assumptions* 2.1, 2.2 *and* 3.1 *hold, then under* $H_{1n}$:

(a) $T_n(r) \Longrightarrow \mu(r) + \sigma G_q(r)$ *in* $D^q[R]$, *and*

(b) $LR_n \overset{\mathcal{L}}{\longrightarrow} \sup_{r \in [a,b]} \{[\sigma^{-1}\mu(r) + G_q(r)]' K_{rr}^{-1} [\sigma^{-1}\mu(r) + G_q(r)]\}$,

*where* $\mu(r) = K_{rr_0} h$ *and* $G_q(r)$ *is a Gaussian process defined as in Theorem* 2.1.

**4. Simulation and one real example.** This section first examines the performance of the statistic $LR_n$ in finite samples through Monte Carlo experiments. In the experiments, sample sizes ($n$) are 200 and 400 and the number of replications is 1000. The null is the MA(1) model with $\phi_{10} = -0.5$ and 0.5 and the alternative is the TMA(1, 1, 2) model with $d = 2$, $r_0 = 0$, $\phi_{10} = 0.5$ and $\psi_{10} = -0.5, -0.3, -0.1, 0.1, 0.3, 0.5$. We take $\beta_1 = 0.1$ and $\beta_2 = 0.9$ in $LR_n$. Significance levels are $\alpha = 0.05$ and 0.1. The corresponding critical values are 7.63 and 9.31, respectively, which were given by Andrews [4]. The results are summarized in Table 1. It shows that the sizes are very close to the nominal values 0.05 and 0.1, in particular, when $n = 400$, and the power increases when the alternative departs from the null model or when the sample size increases. These results indicate that the test has good performance and should be useful in practice.

We next analyze the exchange rate of the Japanese yen against the USA dollar. Monthly data from Jan. 1971 to Dec. 2000 are used and have 360 observations. Define $x_t = 100\Delta \log(\text{exchange rate})$ at the $t$th month and $y_t = x_t - \sum_{t=1}^{360} x_t / 360$. AR(1), TAR(1, 1, 1), MA(1) and TMA(1, 1, 1) models are used to fit the data $\{y_1, \ldots, y_{360}\}$, where the TAR(1, 1, 1) model is defined as in [8]. The results are summarized in Table 2, where $Q(M)$ is the standard Ljung–Box statistic for testing the adequacy of models fitted and $r_0$ is estimated by $\arg\min_{r \in R} \tilde{L}_{1n}(\tilde{\lambda}(r), r)$. The table shows that $Q(11)$, $Q(13)$ and $Q(15)$ all reject AR(1) and TAR(1, 1, 1) models, but they do not reject the MA(1) or TMA(1, 1, 1) models at significance level 0.05.

Based on the MA(1) model, the statistic $LR_n$ is calculated with $\beta_1 = 0.1$ and $\beta_2 = 0.9$ and its value is 14.19. Furthermore, we use the residuals and the estimated $\phi_{10}$ in the MA(1) model to estimate the asymptotic covariance matrix in Theorem 2.2. Using these and the simulation method with 25,000 replications, we obtain that the critical values of the null limiting distribution of $LR_n$ are 6.995, 7.483 and 10.831 at significance levels 0.10, 0.05 and 0.01, respectively. This shows that the null hypothesis of no threshold in the MA(1) model is rejected at all these levels. Furthermore, we note that the

TABLE 1
*Size and power of $LR_n$ for testing threshold in*
*MA(1) models ($\beta_1 = 0.1$, $\beta_2 = 0.9$, $d = 2$,*
*1000 replications)*

|  | $n = 200$ | | $n = 400$ | |
|---|---|---|---|---|
| $\alpha$ | **5%** | **10%** | **5%** | **10%** |
| $\phi_{10}$ | | Sizes | | |
| $-0.5$ | 0.044 | 0.097 | 0.058 | 0.102 |
| 0.5 | 0.059 | 0.112 | 0.051 | 0.101 |
| $\psi_{10}$ | | Powers when $\phi_{10} = 0.5$ | | |
| $-0.5$ | 0.836 | 0.909 | 0.993 | 0.999 |
| $-0.3$ | 0.318 | 0.514 | 0.710 | 0.815 |
| $-0.1$ | 0.076 | 0.156 | 0.123 | 0.191 |
| 0.1 | 0.103 | 0.167 | 0.143 | 0.237 |
| 0.3 | 0.599 | 0.717 | 0.916 | 0.953 |
| 0.5 | 0.989 | 0.993 | 1.000 | 1.000 |

TABLE 2
*Results for monthly exachange rate data of Japanese yen against USA dollar (1971 to*
*2000)*

|  | $\phi_{00}$ | $\psi_{00}$ | $\phi_{10}$ | $\psi_{10}$ | $r_0$ | $Q(11)$ | $Q(13)$ | $Q(15)$ | **AIC** |
|---|---|---|---|---|---|---|---|---|---|
| AR(1) | | | 0.345 | | | 22.66 | 28.91 | 29.26 | 699.83 |
| TAR(1, 1, 1) | 0.930 | $-0.905$ | 0.076 | 0.293 | $-2.51$ | 20.97 | 28.40 | 28.63 | 704.44 |
| MA(1) | | | 0.402 | | | 13.59 | 18.93 | 19.36 | 693.17 |
| TMA(1, 1, 1) | | | 0.281 | 0.445 | $-4.93$ | 15.52 | 19.52 | 19.73 | 691.61 |

Upper-tail 5% critical values: $Q(11) = 19.68$, $Q(13) = 23.36$ and $Q(15) = 25.00$.

TMA$(1,1,1)$ model achieves the minimum AIC among the four candidate models and, hence, it should be a reasonable choice for the data.

Finally, to understand what order of AR or TAR model is adequate for the data, some higher-order models are fitted. We found that AR$(2)$ is not adequate, but AR$(3)$ and TAR$(2,2,1)$ are adequate at significance level 0.05. The result for AR$(3)$ is $y_t = 0.390y_{t-1} - 0.139y_{t-2} + 0.103y_{t-3} + \varepsilon_t$, for which $Q(11) = 13.211$, $Q(13) = 18.106$ and $Q(15) = 18.573$ and the value of AIC is 696.50. The result for TAR$(2,2,1)$ is $y_t = 0.821 + 0.130y_{t-1} - 0.082y_{t-2} + [-0.790 + 0.275y_{t-1} - 0.018y_{t-2}]I(y_{t-1} \le -3.741) + \varepsilon_t$, for which $Q(11) = 12.214$, $Q(13) = 16.936$ and $Q(15) = 17.325$ and the value of AIC is 705.08. In terms of AIC, it is clear that not only are AR$(3)$ and TAR$(2,2,1)$ worse than TMA$(1,1,1)$, they are also worse than MA$(1)$.

**5. Proof of Theorem 2.1.** To prove Theorem 2.1, we first introduce three lemmas. Lemma 5.1 is the basis for the other two lemmas and is similar to Lemma A.1 in [18].

LEMMA 5.1. *If Assumption 2.2 holds, then under $H_0$ it follows that*

(a) $E[|\varepsilon_{t-j}|^k I(r' < y_{t-d} \le r)] \le C(r - r')$ *as* $k = 0, 1, 2, 3, 4$, *and* $j \ge 1$, *and*

(b) $Em_t^k \le C(r - r')$ *as* $k = 1, 2, 3, 4$,

*where* $m_t = \|Z_{t-1}\|I(r' < y_{t-d} \le r)$, $r' < r$, $r, r' \in R_\gamma$, $R_\gamma$ *is defined in Section* 2, *and* $C$ *is a constant independent of* $r'$ *and* $r$.

PROOF. Since $E|\varepsilon_{t-j}|^4 < \infty$, there is a constant $M$ such that $\sup_{|x|>M} |x|^4 \times f(x) < 1$. Since $f$ is continuous, it follows that $\sup_{|x| \le M} |x|^4 f(x) < \infty$. Thus, $\sup_{x \in R} |x|^k f(x) < \infty$ for $k = 0, 1, 2, 3, 4$. Let $g_t = \sum_{i=1}^{p} \phi_{i0}\varepsilon_{t-i}$. When $j \ne d$, $E[|\varepsilon_{t-j}|^k I(r' < y_{t-d} \le r)] = E[|\varepsilon_{t-j}|^k \int_{r'-g_{t-d}}^{r-g_{t-d}} f(x)\,dx] \le C(r - r')$. When $j = d$, $E[|\varepsilon_{t-d}|^k I(r' < y_{t-d} \le r)] = E[\int_{r'-g_{t-d}}^{r-g_{t-d}} |x|^k f(x)\,dx] \le C(r - r')$. Thus, we can show that (a) and (b) hold. $\square$

LEMMA 5.2. *Under the assumptions of Theorem 2.1 and $H_0$, it follows that*:

(a)
$$E\left\| \frac{1}{\sqrt{n}} \sum_{t=1}^{n} \left[ \sum_{i=0}^{\infty} u'\Phi^i u Z_{t-i-1} I(r' < y_{t-d-i} \le r) \right] \varepsilon_t \right\|^4$$
$$\le C\left[ \sqrt{\frac{r - r'}{n}} + (r - r') \right]^2,$$

*and*

(b) $\displaystyle E\left[\frac{1}{\sqrt{n}}\sum_{t=1}^{n}(|\varepsilon_t| - E|\varepsilon_t|)\sum_{i=0}^{\infty}\|\Phi^i\|m_{t-i}\right]^4 \leq C\left[\sqrt{\frac{r-r'}{n}} + (r-r')\right]^2,$

where $C$ is a constant independent of $r'$, $r$ and $n$, and $m_t$ is defined in Lemma 5.1.

PROOF. (a) Let $a_{tj} = \varepsilon_{t-i-j}I(r' < y_{t-d-i} \leq r)$, where $i \geq 0$ and $j = 1, \ldots, q$. Since $\varepsilon_t$ and $a_{tj}$ are independent and $a_{tj}$ is $(p+q+d)$-dependent, we can show that $E(\sum_{t=1}^{n} a_{tj}\varepsilon_t)^4 \leq O(1)\sum_{t=1}^{n}\sum_{t_1=1}^{n} E(a_{tj}^2 a_{t_1 j}^2 \varepsilon_t^2 \varepsilon_{t_1}^2)$, where $O(1)$ holds uniformly in $i$. Note that $\|\Phi^i\| = O(\rho^i)$ with $\rho \in (0,1)$. Thus, by Minkowskii's inequality,

$$
\begin{aligned}
E&\left\|\frac{1}{\sqrt{n}}\sum_{t=1}^{n}\left[\sum_{i=0}^{\infty}u'\Phi^i u Z_{t-i-1}I(r' < y_{t-d-i} \leq r)\right]\varepsilon_t\right\|^4 \\
&= \frac{1}{n^2}E\left\|\sum_{i=0}^{\infty}u'\Phi^i u\sum_{t=1}^{n}[Z_{t-i-1}I(r' < y_{t-d-i} \leq r)]\varepsilon_t\right\|^4 \\
&= \frac{O(1)}{n^2}\left\{\sum_{i=0}^{\infty}\rho^i\left[E\left\|\sum_{t=1}^{n}[Z_{t-i-1}I(r' < y_{t-d-i} \leq r)]\varepsilon_t\right\|^4\right]^{1/4}\right\}^4 \\
&\leq \frac{O(1)}{n^2}\left\{\sum_{i=0}^{\infty}\rho^i\left[E\left(\sum_{t=1}^{n}m_{t-i}^2\varepsilon_t^2\right)^2\right]^{1/4}\right\}^4 \leq O(1)\sum_{i=0}^{\infty}\rho^i E\left(\frac{1}{n}\sum_{t=1}^{n}m_{t-i}^2\varepsilon_t^2\right)^2,
\end{aligned}
$$

$$(5.1)$$

where the third and the last steps hold using the inequality $(\sum_{i=0}^{\infty}\rho^i a_i)^2 = \sum_{i=0}^{\infty}\sum_{j=0}^{\infty}\rho^{i+j}a_i a_j \leq (1-\rho)^{-1}\sum_{i=0}^{\infty}\rho^i a_i^2$, for any $a_i \in R$ as $\sum_{i=0}^{\infty}\rho^i a_i^2 < \infty$. By Lemma 5.1(b),

$$
\begin{aligned}
E\left(\frac{1}{n}\sum_{t=1}^{n}m_{t-i}^2\varepsilon_t^2\right)^2 &= E\left[\frac{1}{n}\sum_{t=1}^{n}m_{t-i}^2(\varepsilon_t^2 - E\varepsilon_{t-i}^2) + E\varepsilon_{t-i}^2\frac{1}{n}\sum_{t=1}^{n}m_{t-i}^2\right]^2 \\
&\leq 2E\left[\frac{1}{n}\sum_{t=1}^{n}m_{t-i}^2(\varepsilon_t^2 - E\varepsilon_{t-i}^2)\right]^2 + 2(E\varepsilon_{t-i}^2)^2 E\left(\frac{1}{n}\sum_{t=1}^{n}m_{t-i}^2\right)^2 \\
&\leq \frac{2}{n^2}\sum_{t=1}^{n}Em_{t-i}^4 E(\varepsilon_t^2 - E\varepsilon_{t-i}^2)^2 + 2(E\varepsilon_{t-i}^2)^2 E\left(\frac{1}{n}\sum_{t=1}^{n}m_{t-i}^2\right)^2 \\
&\leq \frac{C_0(r-r')}{n} + 2(E\varepsilon_{t-i}^2)^2 E\left(\frac{1}{n}\sum_{t=1}^{n}m_{t-i}^2\right)^2,
\end{aligned}
$$

where $C_0$ is a constant independent of $i$, $r'$, $r$ and $n$. Again, by Lemma 5.1(b),

$$E\left(\frac{1}{n}\sum_{t=1}^{n}m_{t-i}^2\right)^2 = E\left[\frac{1}{n}\sum_{t=1}^{n}(m_{t-i}^2 - Em_{t-i}^2) + Em_{t-i}^2\right]^2$$

$$\leq \left\{\frac{1}{n}\left[E\left(\sum_{t=1}^{n}(m_{t-i}^2 - Em_{t-i}^2)\right)^2\right]^{1/2} + C(r-r')\right\}^2.$$

Since $y_t$ is only $p$-dependent, we see that $m_t$ is $\widetilde{p}$-dependent, where $\widetilde{p} = p + q + d$. So, $E[(m_t^2 - Em_t^2)(m_{t_1}^2 - Em_{t_1}^2)] = 0$ when $|t - t_1| > \widetilde{p}$. Thus, by Lemma 5.1(b), it follows that

$$E\left(\sum_{t=1}^{n}(m_{t-i}^2 - Em_{t-i}^2)\right)^2$$

$$= \sum_{t=1}^{n}E(m_{t-i}^2 - Em_{t-i}^2)^2 + 2\sum_{t=1}^{n}\sum_{s=1}^{n-t}E[(m_{t-i}^2 - Em_{t-i}^2)(m_{t-i+s}^2 - Em_{t-i}^2)]$$

$$= \sum_{t=1}^{n}E(m_{t-i}^2 - Em_{t-i}^2)^2$$

$$\quad + 2\sum_{t=1}^{n}\sum_{s=1}^{\min\{n-t,\widetilde{p}\}}E[(m_{t-i}^2 - Em_{t-i}^2)(m_{t-i+s}^2 - Em_{t-i}^2)]$$

$$\leq (2\widetilde{p}+1)\sum_{t=1}^{n}E(m_{t-i}^2 - Em_{t-i}^2)^2 \leq (2\widetilde{p}+1)nC(r-r'),$$

where $C$ is a constant independent of $i$, $r'$, $r$ and $n$. By the preceding three equations and (5.1), we can claim that (a) holds.

(b) Let $\widetilde{\varepsilon}_t = |\varepsilon_t| - E|\varepsilon_t|$. As for (5.1) and the preceeding argument, we have

$$E\left[\frac{1}{\sqrt{n}}\sum_{t=1}^{n}\widetilde{\varepsilon}_t\sum_{i=0}^{\infty}\|\Phi^i\|m_{t-i}\right]^4$$

$$= \frac{O(1)}{n^2}E\left[\sum_{t=1}^{n}\sum_{t_1=1}^{n}\left(\sum_{i=0}^{\infty}\|\Phi^i\|m_{t-i}\right)^2\left(\sum_{i=0}^{\infty}\|\Phi^i\|m_{t_1-i}\right)^2\widetilde{\varepsilon}_t^2\widetilde{\varepsilon}_{t_1}^2\right]$$

$$\leq \frac{O(1)}{n^2}E\left[\sum_{t=1}^{n}\left(\sum_{i=0}^{\infty}\rho^i m_{t-i}\right)^2\widetilde{\varepsilon}_t^2\right]^2 \leq C\left[\frac{(r-r')^{1/2}}{\sqrt{n}} + (r-r')\right]^2,$$

where $C$ is a constant independent of $i$, $r'$, $r$ and $n$. Thus, (b) holds. $\quad\square$

LEMMA 5.3.   *Under the assumptions of Theorem* 2.1 *and* $H_0$, *it follows that*

$$E\left[\frac{1}{\sqrt{n}}\sum_{t=1}^{n}\sum_{i=0}^{\infty}\|\Phi^i\|(m_{t-i}-Em_{t-i})\right]^4 \leq C\left[\frac{r-r'}{n}+(r-r')^2\right],$$

*where* $C$ *is a constant independent of* $r'$, $r$ *and* $n$, *and* $m_t$ *is defined in Lemma* 5.1.

PROOF.   First, for any integer $i \geq 0$, we have the inequality

$$
\begin{aligned}
&E\left[\sum_{t=1}^{n}(m_{t-i}-Em_t)\right]^4 \\
&\quad \leq \sum_{t=1}^{n}E(m_{t-i}-Em_t)^4 + c_1\left|\sum_{t=1}^{n}\sum_{s=1}^{n-t}E[(m_{t-i}-Em_t)^3(m_{t+s-i}-Em_t)]\right| \\
&\quad\quad + c_2\left|\sum_{t=1}^{n}\sum_{s=1}^{n-t}E[(m_{t-i}-Em_t)^2(m_{t+s-i}-Em_t)^2]\right| \\
&\quad\quad + c_3\left|\sum_{t=1}^{n}\sum_{t_1=1}^{n-t}\sum_{t_2=1}^{n-t-t_1}\sum_{t_3=1}^{n-t-t_1-t_2}E[(m_{t-i}-Em_t)\right. \\
&\quad\quad\quad\quad\quad\quad\quad\quad\quad\quad\quad\quad \times(m_{t+t_1-i}-Em_t) \\
&\quad\quad\quad\quad\quad\quad\quad\quad\quad\quad\quad\quad \times(m_{t+t_1+t_2-i}-Em_t) \\
&\quad\quad\quad\quad\quad\quad\quad\quad\quad\quad\quad\quad \left. \times(m_{t+t_1+t_2+t_3-i}-Em_t)]\right| \\
&\quad \equiv A_{1n}+c_1A_{2n}+c_2A_{3n}+c_3A_{4n},
\end{aligned}
$$

(5.2)

*where* $c_1, c_2$ *and* $c_3$ *are constants independent of* $n$ *and* $i$. Since $m_t$ is $\widetilde{p}$-dependent, $E[(m_t-Em_t)^3(m_{t_1}-Em_t)]=0$ when $|t-t_1|>\widetilde{p}$, where $\widetilde{p}=p+q+d$. Thus, by Lemma 5.1(b),

$$
\begin{aligned}
A_{2n} &= \left|\sum_{t=1}^{n}\sum_{s=1}^{\min\{n-t,\widetilde{p}\}}E[(m_{t-i}-Em_t)^3(m_{t+s-i}-Em_{t+s})]\right| \\
&\leq n\widetilde{p}E(m_{t-i}-Em_t)^4 \\
&\leq n\widetilde{p}C_1(r-r').
\end{aligned}
$$

Let $\widetilde{m}_t=(m_{t-i}-Em_t)^2-E(m_{t-i}-Em_t)^2$. Then, by Lemma 5.1(b) we can show that $E\widetilde{m}_t^2\leq C_2(r-r')$. Since $\{\widetilde{m}_t\}$ is a $\widetilde{p}$-dependent sequence, we know that $E(\widetilde{m}_t\widetilde{m}_{t_1})=0$ when $|t-t_1|>\widetilde{p}$. Furthermore, by Lemma 5.1(b),

$$
A_{3n} = \left|\sum_{t=1}^{n}\sum_{s=1}^{n-t}E(\widetilde{m}_t\widetilde{m}_{t+s}) - \sum_{t=1}^{n}(n-t)[E(m_t-Em_t)^2]^2\right|
$$

$$\leq \left| \sum_{t=1}^{n} \sum_{s=1}^{\min\{n-t,\widetilde{p}\}} E(\widetilde{m}_t \widetilde{m}_{t+s}) \right| + C_3 n^2 (r-r')^2$$

$$\leq C_2 \widetilde{p} n (r-r') + C_3 n^2 (r-r')^2.$$

Denote $\widetilde{p}_1 = \min\{n-t, \widetilde{p}\}$. Similarly, by Lemma 5.1(b) we have that

$$A_{4n} = \left| \sum_{t=1}^{n} \sum_{t_1=1}^{\widetilde{p}_1} \sum_{t_2=1}^{\widetilde{p}_1-t_1} \sum_{t_3=1}^{\widetilde{p}_1-t_1-t_2} E[(m_{t-i} - Em_t) \right.$$

$$\times (m_{t+t_1-i} - Em_t)(m_{t+t_1+t_2-i} - Em_t)$$

$$\left. \times (m_{t+t_1+t_2+t_3-i} - Em_t)] \right|$$

$$\leq \widetilde{p}_1^3 \sum_{t=1}^{n} E(m_{t-i} - Em_t)^4$$

$$\leq n \widetilde{p}_1^3 C_4 (r-r').$$

By Lemma 5.1(b), the preceding three inequalities and (5.2), we can claim that

$$E\left[ \sum_{t=1}^{n} (m_{t-i} - Em_t) \right]^4 \leq n C_5 (r-r') + C_5 n^2 (r-r')^2.$$

In the above, $C_i$, $i = 1, \ldots, 5$, are some constants independent of $r'$, $r$, $i$ and $n$. By the assumption given, $\Phi^i = O(\rho^i)$ with $\rho \in (0,1)$. Thus, by Minkowskii's inequality,

$$E\left[ \frac{1}{\sqrt{n}} \sum_{t=1}^{n} \sum_{i=0}^{\infty} \|\Phi^i\| (m_{t-i} - Em_{t-i}) \right]^4$$

$$\leq \frac{1}{n^2} E\left[ \sum_{i=0}^{\infty} \|\Phi^i\| \left| \sum_{t=1}^{n} (m_{t-i} - Em_{t-i}) \right| \right]^4$$

$$\leq \frac{O(1)}{n^2} \left[ \sum_{i=0}^{\infty} \rho^i \left\{ E \left| \sum_{t=1}^{n} (m_{t-i} - Em_{t-i}) \right|^4 \right\}^{1/4} \right]^4$$

$$\leq \frac{O(1)}{n^2} \left\{ [n C_5 (r-r') + C_5 n^2 (r-r')^2]^{1/4} \sum_{i=0}^{\infty} \rho^i \right\}^4$$

$$\leq \frac{C(r-r')}{n} + C(r-r')^2,$$

where $C$ is some constant independent of $r'$, $r$ and $n$. $\quad\square$

PROOF OF THEOREM 2.1. We use Lemmas 5.2 and 5.3 to prove the tightness. Let

$$T_{1n}(r) = \frac{1}{\sqrt{n}} \sum_{t=1}^{n} \left[ \sum_{i=0}^{\infty} u' \Phi^i u Z_{t-i-1} I(y_{t-d-i} \leq r) \right] \varepsilon_t.$$

We first show that $\{T_{1n}(r) : r \in R_\gamma\}$ is tight. For any given $\eta > 0$, we choose $(\delta, n)$ such that $1 > \delta \geq n^{-1}$ and $\sqrt{n} \geq M/\eta$ and then choose an integer $K$ such that $\delta n/2 \leq K \leq n\delta$, where $M$ is determined later.

Let $r_{k+1} = r_k + \delta/K$, where $r_1 = r'$ and $k = 1, \ldots, K$. Thus,

(5.3)
$$\begin{aligned}
\sup_{r' < r \leq r' + \delta} &\|T_{1n}(r) - T_{1n}(r')\| \\
&\leq \sup_{1 \leq k \leq K} \|T_{1n}(r_k) - T_{1n}(r')\| \\
&\quad + \sup_{1 \leq k \leq K} \sup_{r_k < r \leq r_k + \delta/K} \|T_{1n}(r) - T_{1n}(r_k)\|.
\end{aligned}$$

For any $1 \leq i < j \leq K$, we have $(r_j - r_i)^{1/2} = [(j-i)\delta/K]^{1/2} \leq (j-i)\sqrt{\delta/K}$. By Lemma 5.2(a) and the inequality $1/\sqrt{n} \leq \sqrt{\delta/K}$, it follows that

$$E\|T_{1n}(r_i) - T_{1n}(r_j)\|^4 \leq C \left[ \frac{(r_j - r_i)^{1/2}}{\sqrt{n}} + (r_j - r_i) \right]^2 = C \left( \sum_{k=i+1}^{j} \frac{\delta}{K} \right)^2.$$

Note that $T_{1n}(r_j) - T_{1n}(r_i) = \sum_{k=i+1}^{j} [T_{1n}(r_k) - T_{1n}(r_{k-1})]$. By the preceding equation and Theorem 12.2 of [5], page 94, there exists a constant $C_1$ independent of $K$, $\delta$, $r'$ and $n$ such that

(5.4)
$$\begin{aligned}
P\left( \sup_{1 \leq k \leq K} \|T_{1n}(r_k) - T_{1n}(r')\| > \frac{\eta}{2} \right) &\leq \frac{CC_1}{\eta^4} \left( \sum_{k=1}^{K} \frac{\delta}{K} \right)^2 \\
&= \frac{CC_1 \delta^2}{\eta^4}.
\end{aligned}$$

We now consider the second term of the right-hand side in (5.3). Let

$$m_{kt} = \|Z_{t-1}\| I(r_k < y_{t-d} \leq r_k + \delta/K).$$

By Lemma 5.1(b) and the definition of $K$ and $\eta$,

$$\begin{aligned}
\frac{E|\varepsilon_t|}{\sqrt{n}} \sum_{t=1}^{n} E\left( \sum_{i=0}^{\infty} \|\Phi^i\| m_{kt-i} \right) &\leq \frac{C_2 \sqrt{n} \delta}{K} \\
&\leq \frac{2C_2 \sqrt{n} \delta}{n\delta} \\
&= \frac{2C_2}{\sqrt{n}} \leq \frac{\eta}{8},
\end{aligned}$$

as $M \geq 16C_2$, where $C_2$ is a constant independent of $k$, $\delta$, $r'$ and $n$. By the preceding inequality, Lemma 5.3 and Markov's inequality,

$$
\sum_{k=1}^{K} P\left(\frac{1}{\sqrt{n}}\sum_{t=1}^{n}\left[(E|\varepsilon_t|)\sum_{i=0}^{\infty}\|\Phi^i\|m_{kt}\right] > \frac{\eta}{4}\right)
$$
$$
\leq \sum_{k=1}^{K} P\left(\frac{E|\varepsilon_t|}{\sqrt{n}}\sum_{t=1}^{n}\left[\left(\sum_{i=0}^{\infty}\|\Phi^i\|m_{kt}\right) - E\left(\sum_{i=0}^{\infty}\|\Phi^i\|m_{kt}\right)\right] > \frac{\eta}{8}\right)
$$
$$
\leq \frac{C_3}{\eta^4}\sum_{k=1}^{K}E\left[\frac{1}{\sqrt{n}}\sum_{t=1}^{n}\sum_{i=0}^{\infty}\|\Phi^i\|(m_{kt} - Em_{kt})\right]^4
$$
$$
\leq \frac{C_4 K}{\eta^4}\left(\frac{\delta}{nK} + \frac{\delta^2}{K^2}\right)
$$
$$
\leq \frac{2C_4\delta^2}{\eta^4},
$$

since $n^{-1} \leq \delta/K$, where $C_3$ and $C_4$ are constants independent of $K$, $\delta$, $r'$ and $n$. By the preceding inequality, Lemma 5.2(b) and Markov's inequality, we have

(5.5)
$$
P\left(\sup_{1\leq k\leq K}\sup_{r_k<r\leq r_k+\delta/K}\|T_{1n}(r) - T_{1n}(r_k)\| > \frac{\eta}{2}\right)
$$
$$
\leq \sum_{k=1}^{K} P\left(\frac{1}{\sqrt{n}}\sum_{t=1}^{n}\left(|\varepsilon_t|\sum_{i=0}^{\infty}\|\Phi^i\|m_{kt}\right) > \frac{\eta}{2}\right)
$$
$$
\leq \sum_{k=1}^{K} P\left(\frac{1}{\sqrt{n}}\sum_{t=1}^{n}\left[(|\varepsilon_t| - E|\varepsilon_t|)\sum_{i=0}^{\infty}\|\Phi^i\|m_{kt}\right] > \frac{\eta}{4}\right)
$$
$$
+ \sum_{k=1}^{K} P\left(\frac{1}{\sqrt{n}}\sum_{t=1}^{n}\left[(E|\varepsilon_t|)\sum_{i=0}^{\infty}\|\Phi^i\|m_{kt}\right] > \frac{\eta}{4}\right)
$$
$$
\leq \frac{4^4}{\eta^4}\sum_{k=1}^{K}E\left[\frac{1}{\sqrt{n}}\sum_{t=1}^{n}(|\varepsilon_t| - E|\varepsilon_t|)\sum_{i=0}^{\infty}\|\Phi^i\|m_{kt}\right]^4 + \frac{2C_4\delta^2}{\eta^4}
$$
$$
\leq \frac{C_5 K}{\eta^4}\left(\sqrt{\frac{\delta}{nK}} + \frac{\delta}{K}\right)^2 + \frac{2C_4\delta^2}{\eta^4} \leq \frac{(2C_4 + 4C_5)\delta^2}{\eta^4},
$$

since $1/\sqrt{n} \leq \sqrt{\delta/K}$, where $C_5$ is a constant independent of $K$, $\delta$, $r'$ and $n$.

Given $\varepsilon > 0$ and $\eta > 0$, let $\delta = \min\{\varepsilon\eta^4/(2C_4 + 4C_5 + CC_1), 0.5\}$. We first select $M$ such that $M \geq 16C_2$, and then select $N = \max\{\delta^{-1}, M^2/\eta^2\}$. Thus, for any $r' \in R_\gamma$, as $n > N$, by (5.3)–(5.5) it follows that

$$
P\left(\sup_{r'<r\leq r'+\delta}\|T_{1n}(r) - T_{1n}(r')\| > \eta\right) \leq \frac{(2C_4 + 4C_5)\delta^2}{\eta^4} + \frac{CC_1\delta^2}{\eta^4} \leq \delta\varepsilon.
$$

By Theorem 15.5 in [5] (also see the proof of Theorem 16.1 in [5]), we can claim that $\{T_{1n}(r): R_\gamma\}$ is tight. Furthermore, since $\sum_{t=1}^n D_{1t}(\lambda_0, r)/\sqrt{n}$ is tight under $H_0$ and $\Sigma_{1r}$ is continuous in terms of $r$ on $R_\gamma$, we know that $\{T_n(r): R_\gamma\}$ is tight. We can show that the finite-dimensional distributions of $\{T_n(r): r \in R_\gamma\}$ converge weakly to those of $\{\sigma G_q(r): r \in R_\gamma\}$. By Prohorov's theorem in [5], page 37, $T_n(r) \Rightarrow \sigma G_q(r)$ on $D^q[R_\gamma]$ for each $\gamma \in (0, \infty)$. By Theorem 15.5 in [5], almost all the paths of $G_q(r)$ are continuous in terms of $r$. $\square$

**6. Proof of Lemma 2.1.** To prove Lemma 2.1, we need six lemmas. Lemma 6.1 is a basic result. Lemmas 6.3 and 6.4 are for Lemma 2.1(a). Lemmas 6.2 and 6.5 are for Lemma 2.1(b). Lemma 6.6 shows that the effect of initial values is asymptotically ignorable. Most of the results in this section still hold under $H_1$.

LEMMA 6.1. *If Assumption 2.1 holds with $E\varepsilon_t^4 < \infty$, then under $H_0$:*

(a) $E \sup\limits_{\Theta_1} \sup\limits_{r \in [a,b]} \varepsilon_t^4(\lambda, r) < \infty,$

(b) $E \sup\limits_{\Theta_1} \sup\limits_{r \in [a,b]} \left\| \dfrac{\partial \varepsilon_t(\lambda, r)}{\partial \lambda} \right\|^4 < \infty,$

(c) $E \sup\limits_{\Theta_1} \sup\limits_{r \in [a,b]} \left\| \dfrac{\partial^2 \varepsilon_t(\lambda, r)}{\partial \lambda \, \partial \lambda'} \varepsilon_t(\lambda, r) \right\|^2 < \infty.$

PROOF. By Theorem A.2 in the Appendix, under $H_0$ the following expansion holds:

$$(6.1) \quad \varepsilon_t(\lambda, r) = y_t + \sum_{j=1}^\infty u' \prod_{i=1}^j [\Phi + \Psi I(y_{t-d-i+1} \leq r)] u y_{t-j} \qquad \text{a.s.,}$$

where $u$, $\Phi$ and $\Psi$ are defined in Theorem A.2. By (6.1) and Theorem A.1, we have

$$(6.2) \quad \sup_{\Theta_1} \sup_{r \in [a,b]} |\varepsilon_t(\lambda, r)| \leq O(1) \sum_{i=0}^\infty \rho^i |y_{t-i}| \qquad \text{a.s.,}$$

where $\rho \in (0, 1)$. Since $E\varepsilon_t^4 < \infty$, it is readily shown that $Ey_t^4 < \infty$. By Minkowskii's inequality, we can show that $E \sup_{\Theta_1} \sup_{r \in [a,b]} \varepsilon_t^4(\lambda, r) < \infty$. Thus, (a) holds:

$$\frac{\partial \varepsilon_t(\lambda, r)}{\partial \phi_k} = -\varepsilon_{t-k}(\lambda, r) - \sum_{i=1}^p [\phi_i + \psi_i I(y_{t-d} \leq r)] \frac{\partial \varepsilon_{t-i}(\lambda, r)}{\partial \phi_k},$$

$$\frac{\partial \varepsilon_t(\lambda, r)}{\partial \psi_l} = -\varepsilon_{1t-l}(\lambda, r) - \sum_{i=1}^p [\phi_i + \psi_i I(y_{t-d} \leq r)] \frac{\partial \varepsilon_{t-i}(\lambda, r)}{\partial \psi_l},$$

where $\varepsilon_{1t-l}(\lambda, r) = \varepsilon_{t-l}(\lambda, r)I(y_{t-d} \le r)$, $k = 1, \ldots, p$ and $l = 1, \ldots, q$. By Theorem A.2, under $H_0$, the following expansions hold:

$$(6.3) \quad \frac{\partial \varepsilon_t(\lambda, r)}{\partial \phi_k}$$
$$= -\varepsilon_{t-k}(\lambda, r) - \sum_{j=1}^{\infty} u' \prod_{i=1}^{j} [\Phi + \Psi I(y_{t-d-i+1} \le r)] u\varepsilon_{t-k-j}(\lambda, r),$$

$$(6.4) \quad \frac{\partial \varepsilon_t(\lambda, r)}{\partial \psi_l}$$
$$= -\varepsilon_{1t-l}(\lambda, r) - \sum_{j=1}^{\infty} u' \prod_{i=1}^{j} [\Phi + \Psi I(y_{t-d-i+1} \le r)] u\varepsilon_{1t-l-j}(\lambda, r),$$

a.s. Using (6.3) and (6.4), Theorem A.1 and a similar method as for (a), we can show that (b) holds. Similarly, we can show that (c) holds. $\square$

LEMMA 6.2. *If Assumptions* 2.1 *and* 2.2 *hold, then under* $H_0$ $\Omega_r$ *is positive definite for each* $\lambda \in \Theta_1$.

PROOF. It is sufficient to show that if

$$E\left[c'\frac{\partial \varepsilon_t(\lambda, r)}{\partial \lambda}\frac{\partial \varepsilon_t(\lambda, r)}{\partial \lambda'}c\right] = 0,$$

then $c = 0$ for any constant vector $c = (c_1', c_2')'$ with $c_1 = (c_{11}, \ldots, c_{1p})'$ and $c_2 = (c_{21}, \ldots, c_{2q})'$. The above equation holds if and only if $c'\partial \varepsilon_t(\lambda, r)/\partial \lambda = 0$ a.s., from which we can show that

$$\left[\sum_{i=1}^{p} c_{1i}\varepsilon_{t-i}(\lambda, r)\right]I(y_{t-d} > r)$$
$$+ \left[\sum_{i=1}^{p}(c_{1i} + c_{2i})\varepsilon_{t-i}(\lambda, r)\right]I(y_{t-d} \le r) = 0 \qquad \text{a.s.,}$$

where $c_{2i} = 0$ as $i > q$. From this equation, we have that

$$(6.5) \qquad \left[\sum_{i=1}^{p} c_{1i}\varepsilon_{t-i}(\lambda, r)\right]I(y_{t-d} > r) = 0 \qquad \text{a.s.,}$$

$$(6.6) \qquad \left[\sum_{i=1}^{p}(c_{1i} + c_{2i})\varepsilon_{t-i}(\lambda, r)\right]I(y_{t-d} \le r) = 0 \qquad \text{a.s.}$$

Denote the event $A = \{\sum_{i=1}^{p} c_{1i}\varepsilon_{t-i}(\lambda, r) = 0\}$. If $c_{11} \ne 0$, for simplicity let $c_{11} = 1$. Then $A = \{\varepsilon_{t-1}(\lambda, r) = -\sum_{i=2}^{p} c_{1i}\varepsilon_{t-i}(\lambda, r)\}$. Let $g_{1t-1}(\lambda, r) =$

$\sum_{i=1}^p [\phi_i + \psi_i I(y_{t-d-1} \le r)] \varepsilon_{t-i-1}(\lambda, r)$ and $g_{t-2} = g_{1t-1}(\lambda, r) - \sum_{i=1}^p \phi_{i0} \varepsilon_{t-i} - \sum_{i=2}^p c_{1i} \varepsilon_{t-i}(\lambda, r)$:

$$\varepsilon_{t-1}(\lambda, r) = y_t - g_{1t-1}(\lambda, r) = \varepsilon_t + \sum_{i=1}^p \phi_{i0} \varepsilon_{t-i} - g_{1t-1}(\lambda, r)$$

and, hence, $A = \{\varepsilon_{t-1} = g_{t-2}\}$. Since $\varepsilon_{t-1}$ and $g_{t-2}$ are independent and $\varepsilon_t$ has a density function, $P(A) = EI(\varepsilon_{t-1} = g_{t-2}) = E\{E[I(\varepsilon_{t-1} = g_{t-2})|g_{t-2}]\} = 0$. Thus,

$$
\begin{aligned}
P&\left(\left\{\left[\sum_{i=1}^p c_{1i} \varepsilon_{t-i}(\lambda, r)\right] I(y_{t-d} > r) = 0\right\}\right) \\
&= P\left(\left\{\left[\sum_{i=1}^p c_{1i} \varepsilon_{t-i}(\lambda, r)\right] I(y_{t-d} > r) = 0\right\} \cap A^c\right) \\
&= P(\{I(y_{t-d} > r) = 0\} \cap A^c) = P(\{I(y_{t-d} > r) = 0\}) \\
&= P\left(\varepsilon_{t-d} \le r - \sum_{i=1}^p \phi_{i0} \varepsilon_{t-i}\right) = E\left\{\int_{-\infty}^{r - \sum_{i=1}^p \phi_{i0} \varepsilon_{t-i}} f(x)\, dx\right\} > 0,
\end{aligned}
$$

since $f$ is positive, where $f$ is the density of $\varepsilon_t$. This contradicts (6.5). So, $c_{11} = 0$. Similarly, we can show that $c_{12} = \cdots = c_{1p} = 0$. Similarly, we can show that $c_{21} = \cdots = c_{2q}$ using (6.6).  $\square$

LEMMA 6.3. *If Assumptions* 2.1 *and* 2.2 *hold, then under* $H_0$,

$$\inf_{\|\lambda - \lambda_0\| \ge \eta} \inf_{r \in [a,b]} E[\varepsilon_t^2(\lambda, r) - \varepsilon_t^2(\lambda_0, r)] > 0 \qquad \text{for any } \eta > 0.$$

PROOF. Let $V_{t-1}(\lambda, r) = \varepsilon_t(\lambda, r) - \varepsilon_t(\lambda_0, r)$. Then

(6.7)
$$
\begin{aligned}
V_{t-1}(\lambda, r) =& \sum_{i=1}^p [(\phi_i - \phi_{i0}) + (\psi_i - \psi_{i0}) I(y_{t-d} \le r)] \varepsilon_{t-i}(\lambda, r) \\
&+ \sum_{i=1}^p [(\phi_{i0} + \psi_{i0}) I(y_{t-d} \le r)] V_{t-i}(\lambda, r)
\end{aligned}
$$

and, hence, it is independent of $\varepsilon_t$. Note that, under $H_0$, $\varepsilon_t(\lambda_0, r) = \varepsilon_t$. Since $\varepsilon_t(\lambda, r) = \varepsilon_t(\lambda_0, r) + V_{t-1}(\lambda, r)$, we have $E\varepsilon_t^2(\lambda, r) = E\varepsilon_t^2(\lambda_0, r) + EV_{t-1}^2(\lambda, r)$. $EV_{t-1}^2(\lambda, r) = 0$ if and only if $V_{t-1}(\lambda, r) = 0$ a.s. By (6.7) this occurs if and only if $\sum_{i=1}^p [(\phi_i - \phi_{i0}) + (\psi_i - \psi_{i0}) I(y_{t-d} \le r)] \varepsilon_{t-i}(\lambda, r) = 0$ a.s. From the proof of Lemma 6.2, the preceding equation holds if and only if $\lambda = \lambda_0$ for each $r \in [a, b]$. Since $EV_{t-1}^2(\lambda, r)$ is a continuous function of $(\lambda', r)$ and $\Theta_1 \times [a, b]$ is compact, we have $\inf_{\{\|\lambda - \lambda_0\| \ge \eta\} \times [a,b]} EV_{t-1}^2(\lambda, r) > 0$. Thus, the conclusion holds.  $\square$

Lemma 6.4. *If Assumptions* 2.1 *and* 2.2 *hold, then under* $H_0$, *for any* $\varepsilon > 0$,

$$\lim_{n \to \infty} P\left(\frac{1}{n} \sup_{\Theta_1} \sup_{r \in [a,b]} \left| \sum_{t=1}^{n} [\varepsilon_t^2(\lambda, r) - E\varepsilon_t^2(\lambda, r)] \right| > \varepsilon \right) = 0.$$

Proof. Since $\Theta_1$ is compact, we can choose a collection of balls of radius $\delta > 0$ covering $\Theta_1$ and the number of such balls is a finite integer $K_1$. We take a point $\lambda_i$ in the $i$th ball and denote this ball by $V_{\lambda_i}$. Similarly, we divide $[a, b]$ into $K_2$ parts such that $a = r_1 \leq r_2 < \cdots < r_{K_2+1} = b$ with $|r_i - r_{i-1}| \leq \delta$. Thus,

$$P\left(\frac{1}{n} \sup_{\Theta_1} \sup_{r \in [a,b]} \left| \sum_{t=1}^{n} [\varepsilon_t^2(\lambda, r) - E\varepsilon_t^2(\lambda, r)] \right| > \varepsilon \right)$$

$$\leq \sum_{i=1}^{K_1} \sum_{j=1}^{K_2} P\left(\frac{1}{n} \left| \sum_{t=1}^{n} [\varepsilon_t^2(\lambda_i, r_j) - E\varepsilon_t^2(\lambda_i, r_j)] \right| > \frac{\varepsilon}{2} \right)$$

$$+ P\left( \sup_{1 \leq i \leq K_1} \sup_{1 \leq j \leq K_2} \sup_{\lambda \in V_{\lambda_i}} \sup_{r_j < r \leq r_{j+1}} |E[\varepsilon_t^2(\lambda_i, r_j) - \varepsilon_t^2(\lambda, r)]| > \frac{\varepsilon}{4} \right)$$

$$+ P\left(\frac{1}{n} \sup_{1 \leq i \leq K_1} \sup_{1 \leq j \leq K_2} \sup_{\lambda \in V_{\lambda_i}} \sup_{r_j < r \leq r_{j+1}} \left| \sum_{t=1}^{n} [\varepsilon_t^2(\lambda, r) - \varepsilon_t^2(\lambda_i, r_j)] \right| > \frac{\varepsilon}{4} \right)$$

$$\equiv B_{1n} + B_{2n} + B_{3n}, \text{ say.}$$

For any $r' < r$, let $X_t = -\sum_{i=1}^{p} \psi_i I(r' < y_{t-d} \leq r) \varepsilon_{t-i}(\lambda, r')$. By Theorem A.2,

$$\varepsilon_t(\lambda, r) - \varepsilon_t(\lambda, r') = X_t + \sum_{j=1}^{\infty} u' \prod_{i=1}^{j} [\Phi + \Psi I(y_{t-d-i+1} \leq r)] u X_{t-j} \qquad \text{a.s.}$$

By Lemma 5.1(a), we know that $EI(r' < y_{t-d} \leq r' + \delta) = O(\delta)$. Furthermore, by Lemma 6.1(a) and Hölder's inequality, we can show that

$$E \sup_{\lambda \in \Theta_1} \sup_{r' \in [a,b]} \sup_{r' < r \leq r' + \delta} X_t^2 = O(\delta^{1/2}).$$

By the preceding two equations, Theorem A.1 and Minkowskii's inequality, we have

$$E \sup_{\lambda \in \Theta_1} \sup_{1 \leq j \leq K_2} \sup_{r_j < r \leq r_{j+1}} |\varepsilon_t(\lambda, r) - \varepsilon_t(\lambda, r_j)|^2 \leq O(1) \left( \sum_{i=0}^{\infty} \rho^i \delta^{1/4} \right)^2 = O(\delta^{1/2}).$$

By this equation, Lemma 6.1(a) and the Cauchy–Schwarz inequality,

$$(6.8) \qquad E \sup_{\lambda \in \Theta_1} \sup_{1 \leq j \leq K_2} \sup_{r_j < r \leq r_{j+1}} |\varepsilon_t^2(\lambda, r) - \varepsilon_t^2(\lambda, r_j)| = O(\delta^{1/4}).$$

By Taylor's expansion and Lemma 6.1(b), we have

$$E \sup_{1 \le i \le K_1} \sup_{\lambda \in V_{\lambda_i}} \sup_{r \in [a,b]} |\varepsilon_t(\lambda, r) - \varepsilon_t(\lambda_i, r)|^2 \le \delta^2 E \sup_{\Theta_1} \sup_{r \in [a,b]} \left\| \frac{\partial \varepsilon_t(\lambda, r)}{\partial \lambda} \right\|^2 = O(\delta^2).$$

Furthermore, by Lemma 6.1(a) and the Cauchy–Schwarz inequality, we can show that

$$(6.9) \qquad E \sup_{1 \le i \le K_1} \sup_{\lambda \in V_{\lambda_i}} \sup_{r \in [a,b]} |\varepsilon_t^2(\lambda, r) - \varepsilon_t^2(\lambda_i, r)| = O(\delta).$$

By (6.8) and (6.9), we can take $\delta$ small enough such that $B_{2n} = 0$ and

$$B_{3n} \le P\left( \frac{1}{n} \sup_{1 \le i \le K_1} \sup_{\lambda \in V_{\lambda_i}} \sup_{r \in [a,b]} \left| \sum_{t=1}^{n} [\varepsilon_t^2(\lambda, r) - \varepsilon_t^2(\lambda_i, r)] \right| > \frac{\varepsilon}{8} \right)$$

$$+ P\left( \frac{1}{n} \sup_{\lambda \in \Theta_1} \sup_{1 \le j \le K_2} \sup_{r_j < r \le r_{j+1}} \left| \sum_{t=1}^{n} [\varepsilon_t^2(\lambda, r) - \varepsilon_t^2(\lambda, r_j)] \right| > \frac{\varepsilon}{8} \right) < \frac{\varepsilon}{3}.$$

For this $\delta$, $K_1$ and $K_2$ are fixed. By the ergodic theorem, $B_{1n} < \varepsilon/3$ for $n$ large enough. Thus, we can claim that the conclusion holds. $\quad\square$

LEMMA 6.5. *If Assumptions* 2.1 *and* 2.2 *hold, then under* $H_0$, *for any* $\varepsilon > 0$, *there is an* $\eta > 0$ *such that*

$$P\left( \frac{1}{n} \sup_{\|\lambda - \lambda_0\| \le \eta} \sup_{r \in [a,b]} \left\| \sum_{t=1}^{n} [P_t(\lambda, r) - \Omega_r] \right\| > \varepsilon \right) < \varepsilon,$$

*where* $P_t(\lambda, r) = U_t(\lambda, r) U_t'(\lambda, r) + [\partial^2 \varepsilon_t(\lambda, r)/\partial\lambda \, \partial\lambda'] \varepsilon_t(\lambda, r)$.

PROOF. As for Lemma 6.4, the conclusion can be proved by using Lemma 6.1. $\quad\square$

LEMMA 6.6. *If Assumptions* 2.1 *and* 2.2 *hold, then under* $H_0$:

(a) $\dfrac{1}{n} \sup_{\Theta_1} \sup_{r \in [a,b]} \left| \sum_{t=1}^{n} [\varepsilon_t^2(\lambda, r) - \tilde{\varepsilon}_t^2(\lambda, r)] \right| = o_p(1)$,

(b) $\dfrac{1}{\sqrt{n}} \sup_{\Theta_1} \sup_{r \in [a,b]} \left\| \sum_{t=1}^{n} [D_t(\lambda, r) - \tilde{D}_t(\lambda, r)] \right\| = o_p(1)$,

(c) $\dfrac{1}{n} \sup_{\Theta_1} \sup_{r \in [a,b]} \left\| \sum_{t=1}^{n} [P_t(\lambda, r) - \tilde{P}_t(\lambda, r)] \right\| = o_p(1)$,

*where* $P_t(\lambda, r)$ *is defined in Lemma* 6.5 *and typically* $\tilde{D}_t(\lambda, r)$ *is* $D_t(\lambda, r)$ *with the initial values* $y_s = 0$ *for* $s \le 0$.

PROOF. By Lemma 6.1 and Theorem A.1 we can show that the conclusion holds. □

PROOF OF LEMMA 2.1. For any $\eta > 0$, let $c = \inf_{\|\lambda - \lambda_0\| \geq \eta} \inf_{r \in [a,b]} E[\varepsilon_t^2(\lambda, r) - \varepsilon_t^2(\lambda_0, r)]$. By Lemma 6.3 $c > 0$. Furthermore, by Lemma 6.4 we have that

$$P\left(\inf_{r \in [a,b]} \inf_{\|\lambda - \lambda_0\| \geq \eta} \left\{ \sum_{t=1}^n [\varepsilon_t^2(\lambda, r) - \varepsilon_t^2(\lambda_0, r)] - \frac{cn}{2} \right\} < 0 \right)$$

$$= P\left(\inf_{r \in [a,b]} \inf_{\|\lambda - \lambda_0\| \geq \eta} \left\{ \sum_{t=1}^n [\varepsilon_t^2(\lambda, r) - E\varepsilon_t^2(\lambda, r)] \right.\right.$$

$$- \sum_{t=1}^n [\varepsilon_t^2(\lambda_0, r) - E\varepsilon_t^2(\lambda_0, r)]$$

$$\left.\left. + n[E\varepsilon_t^2(\lambda, r) - E\varepsilon_t^2(\lambda_0, r)] - \frac{cn}{2} \right\} < 0 \right)$$

$$\leq P\left(\sup_{r \in [a,b]} \sup_{\Theta_1} \left\{ \left| \frac{1}{n} \sum_{t=1}^n [\varepsilon_t^2(\lambda, r) - E\varepsilon_t^2(\lambda, r)] \right| \right\} > \frac{c}{4} \right) \to 0$$

as $n \to \infty$. Using the preceding equation and Lemma 6.6(a), we can show that

$$P\left(\inf_{r \in [a,b]} \inf_{\|\lambda - \lambda_0\| \geq \eta} \left\{ \sum_{t=1}^n [\tilde{\varepsilon}_t^2(\lambda, r) - \tilde{\varepsilon}_t^2(\lambda_0, r)] - \frac{cn}{4} \right\} < 0 \right) \to 0$$

as $n \to \infty$. Thus, for any $\epsilon > 0$, it follows that

$$P\left(\sup_{r \in [a,b]} \|\tilde{\lambda}_n(r) - \lambda_0\| > \epsilon \right)$$

$$= P\left\{ \|\tilde{\lambda}_n(r) - \lambda_0\| > \epsilon, \sum_{t=1}^n [\tilde{\varepsilon}_t^2(\tilde{\lambda}_n(r), r) - \tilde{\varepsilon}_t^2(\lambda_0, r)] \leq 0, \right.$$

$$\left. \text{for some } r \in [a,b] \right\}$$

$$\leq P\left\{ \inf_{r \in [a,b]} \inf_{\|\lambda - \lambda_0\| > \epsilon} \sum_{t=1}^n [\tilde{\varepsilon}_t^2(\lambda, r) - \tilde{\varepsilon}_t^2(\lambda_0, r)] \leq 0 \right\} \to 0$$

as $n \to \infty$, that is, (a) holds. Using Taylor's expansion, by (a) of this lemma, Lemmas 6.2, 6.5 and 6.6(b)–(c), we can show that (b) holds. For (c), let $D_{1n} = n^{-1/2} \sum_{t=1}^n D_{1t}(\lambda_0, r)$ and $D_{2n} = n^{-1/2} \sum_{t=1}^n D_{2t}(\lambda_0, r)$. $L_{0n}(\tilde{\phi}_n)$ has

the expansion

$$(6.10) \qquad 2[\widetilde{L}_{0n}(\tilde{\phi}_n) - \widetilde{L}_{0n}(\phi_0)] = -D'_{1n}\Sigma^{-1}D_{1n} + o_p(1).$$

By (b) of this lemma and Lemmas 6.5 and 6.6, using Taylor's expansion, it follows that

$$2[\widetilde{L}_{1n}(\tilde{\lambda}_n(r), r) - \widetilde{L}_{1n}(\lambda_0, r)] = -D'_n\Omega_r^{-1}D_n + R_n,$$

where $D_n = [D'_{1n}, D'_{2n}]'$ and $\sup_{r \in [a,b]} |R_n| = o_p(1)$. After some algebra we have

$$2[\widetilde{L}_{1n}(\tilde{\lambda}_n(r), r) - \widetilde{L}_{1n}(\lambda_0, r)] = -T'_n(r)K_{rr}^{-1}T_n(r) - D'_{1n}\Sigma^{-1}D_{1n} + R_n.$$
$$(6.11)$$

Since $\widetilde{L}_{0n}(\phi_0) = \widetilde{L}_{1n}(\lambda_0, r)$ under $H_0$ for each $r$, by (6.10) and (6.11), (c) holds.  $\square$

# APPENDIX

**Invertibility of TMA models.** This appendix gives a general invertible expansion of TMA models, which can be used for TARMA models. We first provide a uniform bound for these coefficients.

THEOREM A.1.   *If Assumption 2.1 holds, then* $\sup_{\Theta_1} \sup_{r \in R} \| \prod_{i=1}^{j}[\Phi + \Psi I(y_{t-i} \leq r)] \| = O(\rho^j)$ *a.s., as* $j \to \infty$*, where* $\rho \in (0,1)$*,*

$$\Phi = \begin{pmatrix} -\phi_1 & \cdots & & -\phi_p \\ I_{p-1} & & O_{(p-1)\times 1} \end{pmatrix} \quad and \quad \Psi = \begin{pmatrix} -\psi_1 & \cdots & & -\psi_p \\ & O_{(p-1)\times p} & \end{pmatrix},$$

*with* $I_k$ *being the* $k \times k$ *identity matrix and* $O_{k \times s}$ *the* $k \times s$ *zero matrix.*

PROOF.   Let $a = \sup_{\Theta_1} \max\{\sum_{i=1}^{p} |\phi_i|, \sum_{i=1}^{p} |\phi_i + \psi_i|\}$. Then  $a \in [0,1]$. Since $\Theta_1$ is compact, if $a = 1$, then there exists a point $\lambda \in \Theta_1$ such that $\sum_{i=1}^{p} |\phi_i| = 1$ or $\sum_{i=1}^{p} |\phi_i + \psi_i| = 1$, which contradicts Assumption 2.1. Thus, $a \in [0,1)$. For any matrix $C = (c_{ij})$, we introduce the notation $|C| = (|c_{ij}|)$. Denote $e_i = (0, \ldots, 0, 1, 0, \ldots, 0)'_{p \times 1}$ with the $i$th element equal to 1, and $v = (1, \ldots, 1)'_{p \times 1}$. Thus,

$$\sup_{r \in R}\left| e_j \prod_{i=1}^{n}[\Phi + \Psi I(y_{t-i} \leq r)]e_k \right|$$

$$\leq \sup_{r \in R} e_j \prod_{i=1}^{n}[|\Phi|I(y_{t-i} > r) + |\Phi + \Psi|I(y_{t-i} \leq r)]v$$

$$\leq \max\left\{ e_j \prod_{i=1}^{n} A_i v : A_i = |\Phi| \text{ or } |\Phi + \Psi| \right\} \qquad \text{a.s.,}$$

for any $j, k = 1, \ldots, p$. It is not difficult to see that $A_n v \leq (a, 1, \ldots, 1)'$, where, for two vectors $B = (b_1, \ldots, b_p)'$ and $C = (c_1, \ldots, c_p)'$, $B \leq C$ means that $b_i \leq c_i$ for $i = 1, \ldots, p$. Since $a \in [0, 1)$, we can see that $A_{n-1} A_n v \leq A_{n-1}(a, 1, \ldots, 1)' \leq (a, a, 1, \ldots, 1)', \ldots,$ and $A_{n-p+1} \ldots A_n v \leq (a, a, \ldots, a)' = av$. Let $n = ps + r$, where $r = 0, 1, \ldots, p - 1$. Then $\sup_{\Theta_1} e_j \prod_{i=1}^n A_i v \leq Ca^s$, where $C > 0$ is a constant independent of $n$. Since $a^s = O[(a^{1/p})^n] = O(\rho^n)$, the conclusion holds. □

THEOREM A.2. *Let* $\{(w_t, y_t) : t \in Z\}$ *be a strictly stationary sequence with* $E|w_t| < \infty$. *If Assumption* 2.1 *holds, then there exists a unique strictly stationary solution* $\{z_t\}$ *to the equation* $z_t = w_t - \sum_{i=1}^p \phi_i z_{t-i} - \sum_{i=1}^q \psi_i I(y_{t-d} \leq r) z_{t-i}$, *with* $p \geq q$, *and* $z_t$ *has the expansion*

$$z_t = w_t + \sum_{j=1}^\infty u' \prod_{i=1}^j [\Phi + \Psi I(y_{t-d-i+1} \leq r)] u w_{t-j},$$

*a.s. and in* $L^1$, *where* $\Phi$ *and* $\Psi$ *are defined as in Theorem* A.1 *and* $u = (1, 0, \ldots, 0)'_{p \times 1}$.

PROOF. Let $\zeta_t = (z_t, \ldots, z_{t-p+1})'$, $A_t = \Phi + \Psi I(y_{t-d} \leq r)$ and $Y_t = u w_t$. We can rewrite $z_t$ in the vector form

(A.1) $$\zeta_t = Y_t + A_t \zeta_{t-1}.$$

We iterate this equation $J$ steps: $\zeta_t = Y_t + \sum_{j=1}^{J-1} \prod_{i=1}^j A_{t-i+1} Y_{t-j} + \prod_{i=1}^J A_{t-i+1} \zeta_{t-J}$. Let $S_J = Y_t + \sum_{j=1}^{J-1} \prod_{i=1}^j A_{t-i+1} Y_{t-j}$. By Theorem A.1 it is not hard to see that

(A.2) $$E\|S_{J_1} - S_{J_2}\| = E\left\| \sum_{j=J_1}^{J_2-1} \prod_{i=1}^j A_{t-i+1} Y_{t-j} \right\|$$
$$\leq O(1) E|w_t| \sum_{j=J_1}^{J_2-1} \rho^j = O(\rho^{J_1})$$

for any $J_1 < J_2$. By (A.2) we can show that $S_J \to S_\infty$ a.s. and in $L^1$. Let $\zeta_t = S_\infty$. Then $\zeta_t$ is a solution of (A.1). To see the uniqueness, suppose that there is another solution $\zeta_t^*$ a.s. and in $L^1$ for model (A.1). Let $V_t = \zeta_t - \zeta_t^*$. $V_t = A_t V_{t-1} = \cdots = \prod_{i=1}^J A_{t-i+1} V_{t-J}$. Since $E\|V_t\| = $ a constant $< \infty$, by Theorem A.1 we can see that $E\|V_t\| = 0$ and, hence, $\zeta_t = \zeta_t^*$ a.s. and in $L^1$. □

**Acknowledgments.** We thank K. S. Chan, J. Pan, two referees, an Associate Editor and the Editor J. Fan for their very helpful comments.

## REFERENCES

[1] AN, H. Z. and CHEN, S. G. (1997). A note on the ergodicity of non-linear autoregressive model. *Statist. Probab. Lett.* **34** 365–372. MR1467442

[2] AN, H. Z. and CHENG, B. (1991). A Kolmogorov–Smirnov type statistic with application to test for nonlinearity in time series. *Internat. Statist. Rev.* **59** 287–307.

[3] AN, H. Z. and HUANG, F. C. (1996). The geometrical ergodicity of nonlinear autoregressive models. *Statist. Sinica* **6** 943–956. MR1422412

[4] ANDREWS, D. W. K. (1993). Tests for parameter instability and structural change with unknown change point. *Econometrica* **61** 821–856. MR1231678

[5] BILLINGSLEY, P. (1968). *Convergence of Probability Measures.* Wiley, New York. MR0233396

[6] BROCKWELL, P., LIU, J. and TWEEDIE, R. L. (1992). On the existence of stationary threshold autoregressive moving-average processes. *J. Time Ser. Anal.* **13** 95–107. MR1165659

[7] CANER, M. and HANSEN, B. E. (2001). Threshold autoregression with a unit root. *Econometrica* **69** 1555–1596. MR1865221

[8] CHAN, K. S. (1990). Testing for threshold autoregression. *Ann. Statist.* **18** 1886–1894. MR1074443

[9] CHAN, K. S. (1991). Percentage points of likelihood ratio tests for threshold autoregression. *J. Roy. Statist. Soc. Ser. B* **53** 691–696. MR1125726

[10] CHAN, K. S. (1993). Consistency and limiting distribution of the least squares estimator of a threshold autoregressive model. *Ann. Statist.* **21** 520–533. MR1212191

[11] CHAN, K. S., PETRUCCELLI, J. D., TONG, H. and WOOLFORD, S. W. (1985). A multiple-threshold AR(1) model. *J. Appl. Probab.* **22** 267–279. MR0789351

[12] CHAN, K. S. and TONG, H. (1985). On the use of the deterministic Lyapunov function for the ergodicity of stochastic difference equations. *Adv. in Appl. Probab.* **17** 666–678. MR0798881

[13] CHAN, K. S. and TONG, H. (1990). On likelihood ratio tests for threshold autoregression. *J. Roy. Statist. Soc. Ser. B* **52** 469–476. MR1086798

[14] CHAN, K. S. and TSAY, R. S. (1998). Limiting properties of the least squares estimator of a continuous threshold autoregressive model. *Biometrika* **85** 413–426. MR1649122

[15] CHEN, R. and TSAY, R. S. (1991). On the ergodicity of TAR(1) processes. *Ann. Appl. Probab.* **1** 613–634. MR1129777

[16] DE GOOIJER, J. G. (1998). On threshold moving-average models. *J. Time Ser. Anal.* **19** 1–18. MR1624163

[17] HANSEN, B. E. (1996). Inference when a nuisance parameter is not identified under the null hypothesis. *Econometrica* **64** 413–430. MR1375740

[18] HANSEN, B. E. (2000). Sample splitting and threshold estimation. *Econometrica* **68** 575–603. MR1769379

[19] KOUL, H. L. (2000). Fitting a two phase linear regression model. *J. Indian Statist. Assoc.* **38** 331–353. MR1923333

[20] KOUL, H. L., QIAN, L. and SURGAILIS, D. (2003). Asymptotics of *M*-estimators in two-phase linear regression models. *Stochastic Process. Appl.* **103** 123–154. MR1947962

[21] KOUL, H. L. and SCHICK, A. (1997). Efficient estimation in nonlinear autoregressive time-series models. *Bernoulli* **3** 247–277. MR1468305

[22] KOUL, H. L. and STUTE, W. (1999). Nonparametric model checks for time series. *Ann. Statist.* **27** 204–236. MR1701108

[23] LING, S. (1999). On the probabilistic properties of a double threshold ARMA conditional heteroskedastic model. *J. Appl. Probab.* **36** 688–705. MR1737046

[24] LING, S. (2003). Weak convergence for marked empirical processes of stationary sequences and goodness-of-fit tests for time series models. Working paper, Dept. Mathematics, Hong Kong Univ. Science and Technology.

[25] LING, S. and MCALEER, M. (2003). On adaptive estimation in nonstationary ARMA models with GARCH errors. *Ann. Statist.* **31** 642–676. MR1983545

[26] LIU, J., LI, W. K. and LI, C. W. (1997). On a threshold autoregression with conditional heteroscedastic variances. *J. Statist. Plann. Inference* **62** 279–300. MR1468167

[27] LIU, J. and SUSKO, E. (1992). On strict stationarity and ergodicity of a nonlinear ARMA model. *J. Appl. Probab.* **29** 363–373. MR1165221

[28] QIAN, L. (1998). On maximum likelihood estimators for a threshold autoregression. *J. Statist. Plann. Inference* **75** 21–46. MR1671666

[29] STUTE, W. (1997). Nonparametric model checks for regression. *Ann. Statist.* **25** 613–641. MR1439316

[30] TONG, H. (1978). On a threshold model. In *Pattern Recognition and Signal Processing* (C. H. Chen, ed.) 101–141. Sijthoff and Noordhoff, Amsterdam.

[31] TONG, H. (1990). *Non-linear Time Series. A Dynamical System Approach.* Clarendon Press, Oxford. MR1079320

[32] TSAY, R. S. (1987). Conditional heteroscedastic time series models. *J. Amer. Statist. Assoc.* **82** 590–604. MR0898364

[33] TSAY, R. S. (1989). Testing and modeling threshold autoregressive processes. *J. Amer. Statist. Assoc.* **84** 231–240. MR0999683

[34] TSAY, R. S. (1998). Testing and modeling multivariate threshold models. *J. Amer. Statist. Assoc.* **93** 1188–1202. MR1649212

[35] WONG, C. S. and LI, W. K. (1997). Testing for threshold autoregression with conditional heteroscedasticity. *Biometrika* **84** 407–418. MR1467056

[36] WONG, C. S. and LI, W. K. (2000). Testing for double threshold autoregressive conditional heteroscedastic model. *Statist. Sinica* **10** 173–189. MR1742107

DEPARTMENT OF MATHEMATICS
HONG KONG UNIVERSITY OF SCIENCE
 AND TECHNOLOGY
CLEAR WATER BAY
HONG KONG
E-MAIL: maling@ust.hk